\documentstyle[11pt,amssymb,amstex]{amsart}
\numberwithin{equation}{section}

\def\bc{{\mathbb C}}

\def\bn{{\mathbb N}}

\def\bq{{\mathbb Q}}

\def\bz{{\mathbb Z}}

\def\g{\gamma}

\def\e{\varepsilon}
\def\l{\lambda}

\def\r{\rho}
\def\s{\sigma}

\newtheorem{thm}{Theorem}[section]
\newtheorem{lem}[thm]{Lemma}
\newtheorem{cor}[thm]{Corollary}
\newtheorem{prop}[thm]{Proposition}

\begin{document}

\title[$p$-adic chaos]
{On chaos of a cubic $p$-adic dynamical system}
\author{Farrukh Mukhamedov}
\address{Farrukh Mukhamedov\\
Departamento de Fisica\\
Universidade de Aveiro\\
Campus Universitario de Santiago\\
3810-193 Aveiro, Portugal}\email{{\tt far75m@@yandex.ru},{\tt
farruh@@fis.ua.pt}}

\author{Jos\'{e} F.F. Mendes}
\address{Jos\'{e} F.F. Mendes\\
Departamento de Fisica \\
Universidade de Aveiro\\
Campus Universit\'{a}rio de Santiago\\
3810-193 Aveiro, Portugal} \email{{\tt jfmendes@@fis.ua.pt}}

\begin{abstract}
In the paper we describe basin of attraction of the $p$-adic
dynamical system $f(x)=x^3+ax^2$. Moreover, we also describe the
Siegel discs of the system, since the structure of the orbits of
the system is related to the geometry of the $p$-adic Siegel
discs. \vskip 0.3cm \noindent {\bf Mathematics Subject
Classification}:
37E99, 37B25, 54H20, 12J12.\\
{\bf Key words}:  Attractor, Siegel disc, $p$-adic dynamics.
\end{abstract}

\maketitle \small
\section{Introduction}

It is known that the $p$-adic numbers were first introduced by the
German mathematician K.Hensel. During a century after their
discovery they were considered mainly objects of pure mathematics.
Starting from 1980's various models described in the language of
$p$-adic analysis have been actively studied. Applications of
$p$-adic numbers in $p$-adic mathematical physics \cite{ADFV, FW,
MP, VVZ}, quantum mechanics \cite{Kh1} and many others \cite{Kh2}
stimulated increasing interest in the study of $p$-adic dynamical
systems. In \cite{Kh2, TVW} $p$-adic field have arisen in physics
in the theory of superstrings, promoting questions about their
dynamics. Also some applications of $p$-adic dynamical systems to
some biological, physical systems were proposed in \cite{ABKO,
AKK1, KhN, WS}. Other studies of non-Archimedean dynamics in the
neighborhood of a periodic and of the counting of periodic points
over global fields using local fields appeared in \cite{CS,HY,L,
P}. Certain rational $p$-adic dynamical systems were investigated
in \cite{KM},\cite{M} which appear from problems of $p$-adic Gibbs
measures \cite{MR1,MR2}. In \cite{B1, B2} the Fatou set of a
rational function defined over some finite extension of $\bq_p$
has been studied. Besides, an analogue of Sullivan's no wandering
domains' theorem for $p$-adic rational functions, which have no
wild recurrent Julia critical points, was proved.

The most studied discrete $p$-adic dynamical systems (iterations
of maps) are the so called monomial systems, i.e. $f(x)=x^n$. In
\cite{AKTS} an asymptotic behavior of such a dynamical system over
$p$-adic and $\bc_p$ was investigated. In \cite{KhN} perturbated
monomial dynamical systems defined by $f_q(x) = x^n+q(x)$, where
the perturbation term $q(x)$ was a polynomial whose coefficients
had small $p$-adic absolute value, was studied. There it was shown
investigated a connection between monomial and perturbated
monomial systems. Formulas for the number of cycles of a specific
length to a given system and the total number of cycles of such
dynamical systems were provided. These investigations show that
the study of perturbated dynamical systems is important. Even for
a quadratic function $f(x)=x^2+c$, $c\in\bq_p$ its chaotic
behavior is complicated (see \cite{Sh,TVW}). In \cite{Sh} the
Fatou and Julia sets of such a $p$-adic dynamical system were
found. Unique ergodicity and ergodicity of monomial and
perturbated dynamical systems have been considered in
\cite{A},\cite{GKL}.

The aim of this paper is to investigate the asymptotic behavior of
a cubic $p$-adic dynamical system $f(x)=x^3+ax^2$ at $|a|_p\neq
1$. Note that globally attracting sets play an important role in
dynamics, restricting the asymptotic behavior to certain regions
of the phase space. However, descriptions of the global attractor
can be difficult as it may contain complicated chaotic dynamics.
Therefore, in the paper we will investigate the basin of
attraction of such a dynamical system. Moreover, we also describe
the Siegel discs of the system, since the structure of the orbits
of the system is related to the geometry of the $p$-adic Siegel
discs (see \cite{AV}).

\section{Preliminaries}

\subsection{$p$-adic numbers}

Let $\bq$ be the field of rational numbers. Throughout the paper
$p$ will be a fixed prime number. Every rational number $x\neq 0$
can be represented in the form $x=p^r\frac{n}{m}$, where
$r,n\in\bz$, $m$ is a positive integer and $p,n,m$ are relatively
prime. The $p$-adic norm of $x$ is given by $|x|_p=p^{-r}$ and
$|0|_p=0$. This norm satisfies so called the strong triangle
inequality
$$
|x+y|_p\leq\max\{|x|_p,|y|_p\}.
$$

From this inequality one can infer that
\begin{eqnarray}\label{inq1}
&&\textrm{if} \ |x|_p\neq |y|_p, \ \ \textrm{then} \ \
|x-y|_p=\max\{|x|_p,|y|_p\}\\[3mm] \label{inq2}&&\textrm{if} \
|x|_p=|y|_p, \ \textrm{then} \ |x-y|_p\leq |2x|_p.
\end{eqnarray}

 This is a ultrametricity of the norm. The completion
of $\bq$ with respect to $p$-adic norm defines the $p$-adic field
which is denoted by $\bq_p$. Note that any $p$-adic number $x\neq
0$ can be uniquely represented in the canonical series:
\begin{equation}\label{rep}
x=p^{\g(x)}(x_0+x_1p+x_2p^2+...) ,
\end{equation}
where $\g=\g(x)\in\bz$ and $x_j$ are integers, $0\leq x_j\leq
p-1$, $x_0>0$, $j=0,1,2,...$ (see more detail \cite{G},\cite{Ko}).
Observe that in this case $|x|_p=p^{-\g(x)}$.

We recall that an integer $a\in \bz$ is called {\it a quadratic
residue modulo $p$} if the equation $x^2\equiv a(\textrm{mod
$p$})$ has a solution $x\in \bz$.

\begin{lem}\label{sq}\cite{VVZ} In order that the equation
\[
x^2=a, \ \ 0\neq a=p^{\g(a)}(a_0+a_1p+...), \ \ 0\leq a_j\leq p-1,
\ a_0>0
\]
has a solution $x\in \bq_p$, it is necessary and sufficient that
the following conditions are satisfied:
\begin{itemize}
\item[(i)] $\g(a)$ is even; \item[(ii)] $a_0$ is a quadratic
residue modulo $p$ if $p\neq 2$, if $p=2$ besides $a_1=a_2=0$.
\end{itemize}
\end{lem}

For any $a\in\bq_p$ and $r>0$ denote
$$
\bar B_r(a)=\{x\in\bq_p : |x-a|_p\leq r\},\ \ B_r(a)=\{x\in\bq_p :
|x-a|_p< r\},
$$
$$
S_r(a)=\{x\in\bq_p : |x-a|_p= r\}.
$$

A function $f:B_r(a)\to\bq_p$ is said to be {\it analytic} if it
can be represented by
$$
f(x)=\sum_{n=0}^{\infty}f_n(x-a)^n, \ \ \ f_n\in \bq_p,
$$ which converges uniformly on the ball $B_r(a)$.

Note the basics of $p$-adic analysis, $p$-adic mathematical
physics are explained in \cite{G,Ko,R,VVZ}

\subsection{Dynamical systems in $\bq_p$}

In this section we recall some known facts about dynamical systems
$(f,B)$ in $\bq_p$, where $f: x\in B\to f(x)\in B$ is an analytic
function and $B=B_r(a)$ or $\bq_p$.

Recall some standard terminology of the theory of dynamical
systems (see for example \cite{PJS}). Let $f:B\to B$ be an
analytic function. Denote $x^{(n)}=f^n(x^{(0)})$, where $x^0\in B$
and $f^n(x)=\underbrace{f\circ\dots\circ f(x)}_n$.
 If $f(x^{(0)})=x^{(0)}$ then $x^{(0)}$
is called a {\it fixed point}. A fixed point $x^{(0)}$ is called
an {\it attractor} if there exists a neighborhood $U(x^{(0)})$ of
$x^{(0)}$ such that for all points $y\in U(x^{(0)})$ it holds
$\lim\limits_{n\to\infty}y^{(n)}=x^{(0)}$, where $y^{(n)}=f^n(y)$.
If $x^{(0)}$ is an attractor then its {\it basin of attraction} is
$$
A(x^{(0)})=\{y\in \bq_p :\ y^{(n)}\to x^{(0)}, \ n\to\infty\}.
$$
A fixed point $x^{(0)}$ is called {\it repeller} if there  exists
a neighborhood $U(x^{(0)})$ of $x^{(0)}$ such that
$|f(x)-x^{(0)}|_p>|x-x^{(0)}|_p$ for $x\in U(x^{(0)})$, $x\neq
x^{(0)}$. For a fixed point  $x^{(0)}$ of a function $f(x)$ a ball
$B_r(x^{(0)})$ (contained in $B$) is said to be a {\it Siegel
disc} if each sphere $S_{\r}(x^{(0)})$, $\r<r$ is an invariant
sphere of $f(x)$, i.e. if $x\in S_{\r}(x^{(0)})$ then all iterated
points $x^{(n)}\in S_{\r}(x^{(0)})$ for all $n=1,2\dots$. The
union of all Siegel discs with the center at $x^{(0)}$ is said to
{\it a maximum Siegel disc} and is denoted by $SI(x^{(0)})$.

{\bf Remark 2.1.} In non-Archimedean geometry, a center of a disc
is nothing but a point which belongs to the disc, therefore, in
principle, different fixed points may have the same Siegel disc.

Let $x^{(0)}$ be a fixed point of an analytic function $f(x)$. Set
$$
\l=\frac{d}{dx}f(x^{(0)}).
$$

The point $x^{(0)}$ is called {\it attractive} if $0\leq |\l|_p<1$, {\it indifferent} if
$|\l|_p=1$, and {\it repelling} if $|\l|_p>1$.

\section{The map $f:x\to x^3+ax^2$ and its fixed points}

In this section we consider dynamical system associated with the
function $f:\bq_p\to\bq_p$ defined by
\begin{equation}\label{func}
f(x)=x^3+ax^{2}, \ \ \ \ a\in\bq_p.
\end{equation}

Throughout the paper we will consider only the case $|a|_p\neq 1$,
since the case $|a|_p=1$ requires more investigations and it will
be considered elsewhere.

Direct checking shows that the fixed points of the function
\eqref{func} are the following ones
\begin{equation}\label{fix}
x_{1}=0 \ \ \textrm{and} \ \ x_{2,3}=\frac{-a\pm\sqrt{a^2+4}}{2}.
\end{equation}

Here it should be noted that $x_{2,3}$ are the solutions of the
following equation
\begin{equation}\label{eq1}
x^2+ax-1=0.
\end{equation}

Note that these fixed points are formal, because, basically in
$\bq_q$ the square root does not always exist, but a full
investigation of behavior of the dynamics of the function needs
the existence of the fixed points $x_{2,3}$. To do end this it is
enough to verify when $\sqrt{a^2+4}$ does exist. So, by verifying
the conditions of Lemma \ref{sq} one can prove the following

\begin{prop}\label{11} The following assertions hold
    \begin{itemize}
    \item[(i)] Let $|a|_p<1$, then the expression $\sqrt{a^2+4}$ exists in
$\bq_p$ if and only if either $p\geq 3$ or $p=2$ and $|a|_p\leq
1/p^3$.
    \item[(ii)] Let $|a|_p>1$, then the expression $\sqrt{a^2+4}$ exists in $\bq_p$.
    \end{itemize}
\end{prop}

\section{Attractors and Siegel discs}

In this  section we will establish behavior of the fixed points of
the dynamical system. After, we are going to describe size of the
attractors and Siegel discs of the system.

Before going to details let us formulate certain useful results.

Let us assume that $x^{(0)}$ be a fixed point of $f$. Then $f$ can
be represented as follows
\begin{equation}\label{tay}
f(x)=f(x^{(0)})+f'(x^{(0)})(x-x^{(0)})+\frac{f''(x^{(0)})}{2}(x-x^{(0)})^2+\frac{f'''(x^{(0)})}{6}(x-x^{(0)})^3.
\end{equation}
From the above equality putting $\g=x-x_0$ we obtain
\begin{equation}\label{dif1}
|f(x)-f(x^{(0)})|_p=|\g|_p\bigg|f'(x^{(0)})+\frac{f''(x^{(0)})}{2}\g+\frac{f'''(x^{(0)})}{6}\g^2\bigg|_p.
\end{equation}

\begin{lem}\label{att1}
Let $x^{(0)}$ be a fixed point of the function $f$ given by
\eqref{func}. If
\begin{equation*}
\max\bigg\{|3x^{(0)}+a|_p|x-x^{(0)}|_p,|x-x^{(0)}|^2_p\bigg\}<|f'(x^{(0)})|_p
\end{equation*}
then
\begin{equation}\label{dif2}
|f(x)-f(x^{(0)})|_p=|f'(x^{(0)})|_p|x-x^{(0)}|_p.
\end{equation}
\end{lem}

The proof immediately follows from \eqref{dif1} and
\begin{equation*} f''(x)=6x+2a, \ \ \ \ f'''(x)=6.
\end{equation*}
From Lemma \ref{att1} we get

\begin{cor}\label{att2}\cite{AKTS} Let $x^{(0)}$ be a fixed point of the function $f$ given by
\eqref{func}. The following assertions hold:
\begin{itemize}
    \item[(i)] if $x^{(0)}$ is an attractive point of $f$, then it is an attractor
of the dynamical system. If $r>0$ satisfies the inequality
\begin{equation}\label{att3}
\max\bigg\{|3x^{(0)}+a|_pr,r^2\bigg\}<1
\end{equation}
then $B_r(x^{(0)})\subset A(x^{(0)})$;
    \item[(ii)] if $x^{(0)}$ is an indifferent point of $f$ then it is the
center of a Siegel disc. If $r$ satisfies the inequality
\eqref{att3} then $B_r(x^{(0)})\subset SI(x^{(0)})$;
    \item[(iii)] if $x^{(0)}$ is a repelling point of $f$ then $x^{(0)}$ is a
repeller of the dynamical system.
\end{itemize}
\end{cor}

Now we going to calculate norms of the fixed points and their
behavior. Consider several distinct cases with respect to the
parameter $a$.

Let us first note that the derivative of $f$ is
\begin{equation}\label{der1}
f'(x)=3x^2+2ax.
\end{equation}

From this one concludes that the fixed point $x_1$ is attractive.
Therefore, furthermore we will deal with $x_{2,3}$.

Taking into account that the fixed points $x_2$ and $x_3$ are the
solutions of \eqref{eq1} from \eqref{der1} we find
\begin{equation}\label{der2}
f'(x_\s)=3-ax_\s, \qquad \s=2,3,
\end{equation}
and
\begin{equation}\label{viet}
x_2+x_3=-a, \ \ \ \ x_2x_3=-1.
\end{equation}

{\tt Case $|a|_p>1$}

In this case from \eqref{viet} one gets that $|x_2+x_3|_p=|a|_p$,
$|x_2|_p|x_3|_p=1$. These imply that either $|x_2|_p>1$ or
$|x_3|_p>1$. Without loss of generality we may assume that
$|x_2|_p>1$, which means that $|x_3|_p<1$.  From the condition
$|a|_p>1$ one finds that $|x_2|_p=|a|_p$ and $|x_3|_p=1/|a|_p$.

From \eqref{der2} we infer that
\begin{equation}\label{der6}
|f'(x_2)|_p=|a|_p|x_2|_p=|a|_p^2>1
\end{equation}
which implies that the point $x_2$ is repelling.

From \eqref{der1} with $|x_3|_p=1/|a|_p$ one gets that
$|f'(x_3)|_p=1$, hence $x_3$ is a indifferent point.

Now let us find the basin of attraction of $x_1$, i.e. $A(x_1)$.

Denote
$$
r_k=\frac{1}{|a|^k_p}, \ \ k\geq 0.
$$

Now consider several steps along the description of $A(x_1)$.

{\tt (I)}. From Corollary \ref{att2} and \eqref{att3} we find that
$B_{r_1}(0)\subset A(x_1)$. Now take $x\in S_{r_1}(0)$, i.e.
$|x|=r_1$. Then one gets
$$
|f(x)|_p=|x|^2_p|x+a|_p=|x|^2_p|a|_p=r_1,
$$
whence we infer that $|f^{(n)}(x)|_p=r_1$ for all $n\in\bn$. This
means that $x\not\in A(x_1)$, hence $A(x_1)\cap
S_{r_1}(0)=\emptyset$. As a consequence we have
$f(S_{r_1}(0))\subset S_{r_1}(0)$.

In the sequel we will assume that $\sqrt{|a|_p}\notin\{p^k, \
k\in\bn\}.$ Denote
$$
A(\infty)=\{x\in\bq_p: |f^{(n)}(x)|_p\to\infty \ \ \textrm{as} \ \
n\to\infty\}.
$$
It is evident that $A(x_1)\cap A(\infty)=\emptyset$.

{\tt (II)}. Let us take $x\in S_{r}(0)$ with $r>|a|_p$. Then we
have
$$
|f(x)|_p=|x|^2_p|x+a|_p=|x|^2_p|x|_p=|x|^3_p,
$$
which means that $x\in A(\infty)$, i.e. $S_{r}(0)\subset
A(\infty)$ for all $r>|a|_p$.

{\tt (III)}. Now assume that $x\in S_r(0)$ with $r\in
(r_1,r_0)\cup (r_0,|a|_p)$. Then we have $f(S_r(0))\subset
S_{r^2|a|_p}(0)$. If $r\in (r_0,|a|_p)$ then $r^2|a|_p>|a|_p$,
hence we have $S_r(0)\subset A(\infty)$. If $r\in (r_1,r_0)$ then
according to our assumption we have
$r^{2^n}|a|_p^{1+2+\cdots+2^{n-1}}\neq 1$ for every $n\in\bn$,
hence there is $n_0\in\bn$ such that $f^{(n_0)}(S_r(0))\subset
A(\infty)$, from this we infer that $S_r(0)\subset A(\infty)$.
Consequently, we have $S_{r}(0)\subset A(\infty)$ for all $r\in
(r_1,r_0)\cup (r_0,|a|_p)$.

{\tt (IV)}. If $x\in S_{r_0}(0)$, then one gets
$f(S_{r_0}(0))\subset S_{|a|_p}(0)$.

{\tt (V)}. Therefore, we have to consider $x\in S_{|a|_p}(0)$.
From \eqref{func} we can write
\begin{equation}\label{ff}
|f(x)|_p=|a|^2_p|x+a|_p.
\end{equation}

From the last equality and the following decomposition
\begin{equation}\label{dec}
S_{|a|_p}(0)=\bigcup_{r=0}^{|a|_p}S_{r}(-a).
\end{equation}
one concludes that we have to investigate behavior of $f$ on
spheres $S_r(-a)$ ($r\in[0,|a|_p]$).

{\tt (VI)}. Let $x\in B_{r_3}(-a)$, i.e. $|x+a|_p<r_3$, then the
equality \eqref{ff} implies that $|f(x)|_p<r_1$ which means $x\in
A(x_1)$. Hence $B_{r_3}(-a)\subset A(x_1)$. Moreover, taking into
account (I) one gets $f(S_{r_3}(-a))\subset S_{r_1}(0)$.

{\tt (VII)}. If $x\in S_{r_1}(-a)$ then from \eqref{ff} we find
that $f(x)\in S_{|a|_p}(0)$.

{\tt (VIII)} If $x\in S_{r_2}(-a)$ then again using \eqref{ff} one
gets that $f(x)\in S_{r_0}(0)$. This with (IV) implies that
$f^{(2)}(S_{r_1}(-a))\subset S_{|a|_p}(0)$.

{\tt (IX)} If $x\in S_{r}(-a)$ with $r\in (r_3,r_2)\cup
(r_2,r_1)\cup (r_1,|a|_p]$. Then from \eqref{ff} we obtain that
$f(x)\in S_{\rho}(0)$, $\rho\in (r_1,r_0)\cup (r_0,|a|_p)\cup
(|a|_p,|a|^3_p]$. Hence thanks to (II) and (III) we infer that
$f(x)\in A(\infty)$.

Let us introduce some more notations. Given sets $A,B\subset\bq_p$
put
\begin{eqnarray}\label{TD1}
T_{f,A,B}(x)&=&\min\{k\in\bn: \ f^{(k)}(x)\in B\}, \ \ x\in A, \\
\label{TD2} D[A,B]&=&\{x\in A: \ T_{f,A,B}(x)<\infty\}.
\end{eqnarray}

Taking into account (II)-(IX) and \eqref{TD1}-\eqref{TD2} we can
define $D[S_{r_0}(0)\cup S_{|a|_p}(0),B_{r_3}(-a)]$, which is non
empty, since from (VI) one sees that $B_{r_3}(-a)\subset
D[S_{r_0}(0)\cup S_{|a|_p}(0),B_{r_3}(-a)]$. Thus from (I)-(IX) we
have $A(x_1)=B_{r_1}(0)\cup D[S_{r_0}(0)\cup
S_{|a|_p}(0),B_{r_3}(-a)]$.

Now consider the other $x_\s$ ($\s=2,3)$ fixed points. We already
have known from above calculations that $|x_2|_p=|a|_p$
and$|x_3|_p=r_1$, this means $x_3$ is indifferent, so Corollary
\ref{att2} with \eqref{att3} yields that $B_{r_1}(x_3)\subset
SI(x_3)$. It is clear that $0\notin SI(x_3)$, therefore
$SI(x_3)=B_{r_1}(x_3)$. From (I) we infer that $SI(x_3)\subset
S_{r_1}(0)$.

Thus we have proved the following

\begin{thm} Let $|a|_p>1$ and $\sqrt{|a|_p}\notin\{p^k, \
k\in\bn\}.$. Then $x_1$ is attractor and $A(x_1)=B_{r_1}(0)\cup
D[S_{r_0}(0)\cup S_{|a|_p}(0),B_{r_3}(-a)]$. For the other fixed
points we have $|x_2|_p=|a|_p$ and $|x_3|_p=r_1$, hence  $x_3$ is
indifferent and $SI(x_3)=B_{r_1}(x_3)$.
\end{thm}

{\tt Case $|a|_p<1$}.

Let $p\geq 3$, from \eqref{fix} and \eqref{viet} one finds that
$|x_\s|_p=1$.

Let $p=2$, then Proposition \ref{11} implies that $a=p^k\e$ for
some $k\geq 3$ and $|\e|_p=1$. From this and taking into account
\eqref{fix} we have
\begin{eqnarray*}
|x_\s|_p=|p^{k-1}\e\pm\sqrt{p^{2(k-1)}\e^2+1}|_p=1,
\end{eqnarray*}
since $|p^{2(k-1)}\e^2+1|_p=1$ and $k\geq 3$.

Now let us compute $|f'(x_\s)|_p$.  If $p\neq 3$, then using
\eqref{der2} one gets
\begin{equation*}
|f'(x_\s)|_p=|3-ax_\s|_p=1, \ \ \s=2,3.
\end{equation*}
This means that the fixed points $x_{\s}$,($\s=2,3$) are
indifferent.

If $p=3$, then from \eqref{der2} we easily obtain that
$|f'(x_\s)|_p<1$, $ \s=2,3$, which implies that the fixed points
are attractive.

Let us first consider $x_1$. According to Corollary \ref{att2} we
immediately find that $B_1(0)\subset A(x_1)$. Take $x\in S_1(0)$
then $|f(x)|_p=|x|^2_p|x+a|_p=|x|^3=1$, hence $|f^{(n)}(x)|_p=1$
for all $n\in\bn$. This means that $x\not\in A(x_1)$, hence
$A(x_1)=B_{1}(0)$.

In the sequel according to the above done calculations we will
consider two possible situations when $p\neq 3$ and $p=3$.

Assume $p\neq 3$. In this case $x_2$ and $x_3$ are indifferent, so
Corollary \ref{att2} implies that $B_{1}(x_\s)\subset SI(x_\s)$,
$\s=2,3$.

Let us take $x\in S_r(x_\s)$, $r\geq 1$, then put $\g=x-x_\s$. It
is clear that $|\g|_p=r$. By means of \eqref{dif1} and \eqref{eq1}
we find
\begin{eqnarray}\label{dif3}
|f(x)-f(x_\s)|_p&=&|\g|_p|3x_\s^2+2ax_\s+(3x_\s+a)\g+\g^2|_p
\nonumber
\\
&=&r|\g^2+3x_\s\g+3+a(\g-x_\s)|_p.
\end{eqnarray}

If $r>1$ then from we easily obtain that
\[
|f(x)-f(x_\s)|_p=|\g|_p^3
\]
since $|\g^2+3x_\s\g+3|_p=|\g|_p^2$, $|a(\g-x_\s)|_p=|a|_p|\g|_p$.
This implies that $SI(x_\s)\subset \bar B_1(x_\s)$.

\begin{lem}\label{sig}
Let $|a|_p<1$ and $p\neq 3$. The equality $SI(x_\s)=\bar
B_1(x_\s)$ holds if and only if for every $\g\in S_1(0)$ the
equality
\begin{equation}\label{sig1}
|\g^2+3x_\s\g+3|=1
\end{equation}
is valid.
\end{lem}

\begin{pf} If \eqref{sig1} is satisfied for all
$\g\in S_1(0)$ then from \eqref{dif3} we infer that
$f(S_1(x_\s))\subset S_1(x_\s)$, since $|a(\g-x_\s)|_p<1$. This
proves the assertion. Now suppose that $SI(x_\s)=\bar B_1(x_\s)$
holds. Assume that \eqref{sig1} is not valid, i.e. there is
$\g_0\in S_1(0)$ such that
\begin{equation}\label{sig2}
|\g_0^2+3x_\s\g_0+3|<1
\end{equation}
The last one with \eqref{dif3} implies that $|f(x_0)-x_\s|_p<1$
for an element $x_0=x_\s+\g_0$. But this contradicts to
$SI(x_\s)=\bar B_1(x_\s)$. This completes the proof.
\end{pf}

From the proof of Lemma \ref{sig} we immediately obtain that if
there is $\g_0\in S_1(0)$ such that \eqref{sig2} is satisfied then
$SI(x_\s)=B_1(x_\s)$.

\begin{lem}\label{sig22} Let $|a|_p<1$ and $p\neq 3$. The following conditions are
equivalent:
\begin{itemize}
    \item[(i)] $SI(x_\s)=B_1(x_\s)$;
    \item[(ii)] there is $\g_0\in S_1(0)$ such that \eqref{sig2} is
    satisfied;
    \item[(iii)] $\sqrt{-3}$ exists in $\bq_p$.
\end{itemize}
\end{lem}

\begin{pf} The implication (i)$\Leftrightarrow$(ii) immediately
follows from the proof of Lemma \ref{sig1}. Consider the
implication (ii)$\Rightarrow$(iii). The condition \eqref{sig2}
according to the Hensel Lemma yields that existence of the
solution $z\in\bq_p$ of the following equation
\begin{equation}\label{sig3}
z^2+3x_\s z+3=0
\end{equation}
such that $|z-\g_0|<1$ which implies that $|z|_p=1$. Now assume
that there is a solution $z_1\in\bq_p$ of \eqref{sig3}. Then from
Vieta's formula we infer the existence of the other solution
$z_2\in\bq_p$ such that
$$
z_1+z_2=-3x_\s, \ \ \ z_1z_2=3.
$$
From these equalities we obtain that $|z_1+z_2|_p=1$,
$|z_+z_-|_p=1$ which imply $z_1,z_2\in S_1(0)$. So putting
$\g_0=z_1$ one gets \eqref{sig2}.

Let us now analyze when \eqref{sig3} has a solution belonging to
$\bq_p$. We know that a general solution of \eqref{sig3} is given
by
\begin{equation}\label{sig4}
z_{1,2}=\frac{-3x_\s\pm\sqrt{-3-9ax_\s}}{2},
\end{equation}
here we have used \eqref{eq1}. But it belongs to $\bq_p$ if
$\sqrt{-3-9ax_\s}$ exists in $\bq_p$. Since $|-9ax_\s|_p=|a|_p<1$
implies that $-9ax_\s=p^k\e$ for some $k\geq 1$ and $|\e|_p=1$.
Hence, $-3-9ax_\s=-3+p^k\e$. Therefore according to Lemma \ref{sq}
we conclude that $\sqrt{-3-9ax_\s}$ exists if and only if
$\sqrt{-3}$ exists in $\bq_p$. The implication
(iii)$\Rightarrow$(ii) can be proven along the reverse direction
in the previous implication.
\end{pf}

From Lemmas \ref{sig} and \ref{sig22} we conclude that $SI(x_\s)$
is either $B_1(x_\s)$ or $\bar B_1(x_\s)$. The equality
\eqref{fix} yields that
\begin{equation}\label{sol}
|x_2-x_3|_p=|\sqrt{a^2+4}|_p=|2|_p,
\end{equation}
which implies that $SI(x_2)\cap SI(x_3)=\emptyset$ when $p\geq 5$
and $SI(x_2)=SI(x_3)$ when $p=2$, since any point of a ball is its
center in non-archimedean setting.\\

Now consider the case $p=3$.  Then, we know that both fixed points
$x_2$ and $x_3$ are attractive. Taking into account $|x_\s|_p<1$
and $|a|_p<1$ from Corollary \ref{att2} one finds that
$B_{1}(x_\s)\subset A(x_\s)$, $\s=2,3$. From the equality
\eqref{sol} we have $|x_2-x_3|_p=1$ which implies that
$S_1(x_\s)\nsubseteq A(x_\s)$.

Let us take $x\in S_r(x_\s)$ with $r\geq 1$, then putting
$\g=x-x_\s$ from \eqref{dif3} with
$|3x_\s\g+3+a(\g-x_\s)|_p<|\g|_p$ we get
\begin{eqnarray*}
|f(x)-x_\s|_p&=&|\g||\g^2+3x_\s\g+3+a(\g-x_\s)|_p=r^3,
\end{eqnarray*}
which implies that $f(S_r(x_\s))\subset S_{r^3}(x_\s)$ for every
$r\geq 1$. Hence, in particular, we obtain $f(S_1(x_\s))\subset
S_{1}(x_\s)$.

Consequently we have the following

\begin{thm}\label{sig5} Let $|a|_p<1$. The following assertions hold:
\begin{itemize}
    \item[(i)] The fixed point $x_1$ is attractor and
    $A(x_1)=B_1(0)$;
    \item[(ii)] If $p\neq 3$ the fixed points $x_\s$,$\s=2,3$ are indifferent
    and $SI(x_\s)=B_1(x_\s)$ is valid if and only if $\sqrt{-3}$
    exists in $\bq_p$. Otherwise $SI(x_\s)=\bar B_1(x_\s)$ holds.
    \item[(iii)] If $p\geq 5$ then $SI(x_2)\cap SI(x_3)=\emptyset$, if
    $p=2$ then $SI(x_2)=SI(x_3)$.
    \item[(iv)] If $p=3$ then the fixed points $x_\s$, $\s=2,3$ are
    attractors and $A(x_\s)=B_1(x_\s)$.
\end{itemize}
\end{thm}

Note that if we consider our dynamical system over $p$-adic
complex field $\bc_p$ we will obtain different result from the
formulated Theorem, since $\sqrt{-3}$ always exists in $\bc_p$.


\section*{ Acknowledgements}  F.M. thanks the FCT (Portugal) grant SFRH/BPD/17419/2004.
J.F.F.M. acknowledges projects POCTI/FAT/46241/2002,
POCTI/MAT/46176/2002 and European research NEST project DYSONET/
012911.

\end{document}